\documentclass[letterpaper]{article} 
\usepackage{aaai23}  
\usepackage{times}  
\usepackage{helvet}  
\usepackage{courier}  
\usepackage[hyphens]{url}  
\usepackage{graphicx} 
\usepackage{natbib}  
\usepackage{caption} 
\usepackage{subcaption}
\usepackage{paralist}
\usepackage{url}
\usepackage{verbatim}

\usepackage{xcolor}

\definecolor{darkgrn}{rgb}{0, 0.75, 0}

\usepackage{xargs}                      
\usepackage[colorinlistoftodos,prependcaption,textsize=footnotesize]{todonotes}
\usepackage{xspace}

\newcommandx{\ak}[2][1=]{\todo[color=green!50,#1]{\sf \textbf{AK:} #2}\xspace}

\newcommandx{\tz}[2][1=]{\todo[color=green!50,#1]{\sf \textbf{BK:} #2}\xspace}

\newcommandx{\ob}[2][1=]{\todo[color=green!50,#1]{\sf \textbf{PP:} #2}\xspace}

\newcommandx{\pg}[2][1=]{\todo[color=green!50,#1]{\sf \textbf{SW:} #2}\xspace}

\usepackage{algorithm}
\usepackage{algorithmic}

\usepackage{newfloat}
\usepackage{listings}
\DeclareCaptionStyle{ruled}{labelfont=normalfont,labelsep=colon,strut=off} 
\floatstyle{ruled}
\newfloat{listing}{tb}{lst}{}
\floatname{listing}{Listing}
\usepackage{makecell}
\usepackage{chngcntr,tocloft}

\title{Persistent Homology to Study Cold Hardiness of Grape Cultivars}
\author{
    Sejal Welankar,\textsuperscript{\rm 1}
    Paola Pesantez-Cabrera,\textsuperscript{\rm 1}
    Bala Krishnamoorthy,\textsuperscript{\rm 2}
    Lynn Mills,\textsuperscript{\rm 3}
    Markus Keller,\textsuperscript{\rm 3}
    Ananth Kalyanaraman\textsuperscript{\rm 1}
}
\affiliations{
    \textsuperscript{\rm 1}School of Electrical Engineering and Computer Science, Washington State University, Pullman, WA 99164\\
    \textsuperscript{\rm 2}Department of Mathematics and Statistics, Washington State University, Vancouver, WA 98686\\
    \textsuperscript{\rm 3}Irrigated Agriculture Research and Extension
Center, Washington State University, Prosser, WA 99350\\

    sejal.welankar@wsu.edu, p.pesantezcabrera@wsu.edu, kbala@wsu.edu, ljmills@wsu.edu, mkeller@wsu.edu, ananth@wsu.edu
}

\usepackage{bibentry}

\begin{document}

\maketitle

\begin{abstract}
Persistent homology is a branch of computational algebraic topology that studies shapes and extracts features over multiple scales. In this paper, we present an unsupervised approach that uses persistent homology to study divergent behavior in agricultural point cloud data. More specifically, we build persistence diagrams from multidimensional point clouds, and use those diagrams as the basis to compare and contrast different subgroups of the population. We apply the framework to study the cold hardiness behavior of 5 leading grape cultivars, with real data from over 20 growing seasons. Our results demonstrate that persistent homology is able to effectively elucidate divergent behavior among the different cultivars; identify cultivars that exhibit variable behavior across seasons; and identify seasonal correlations.

\end{abstract}

\section{Introduction}
\label{sec:intro}

Multidimensional point cloud data sets are becoming pervasive in numerous agricultural applications. Advances in phenotyping technologies that measure various crop attributes, coupled with an increased adoption of field sensing for environmental monitoring (e.g., temperature, humidity, soil moisture) have led to increased availability of multidimensional data sets in agriculture.
While most of these variables are temporal, some variables may also encode static attributes that describe spatial locators or the crop varieties (cultivars) grown at those locations.

Given such a complex spatiotemporal data set, we are interested in understanding how different cultivars (with different genotypes G) respond to different environmental factors (E) to affect their phenotypes (P) \cite{tardieu_plant_2017}. This question of decoding the G$\times$E$\rightarrow$ P interaction is at the center of modern-day phenomics.
However, prior to understanding this complex interaction, historical data sets offer a more immediate opportunity to understand how different genotypes relate to one another by their phenotypic behavior. For instance, extracting different patterns of phenotypic behavior---be it conserved or divergent---among various subgroups of cultivars could help classify cultivars into  behavioral groups. Class information could help in subsequent prediction tasks associated with those cultivars. However, in agricultural data sets, such a classification may not be trivially observable from data, particularly for temporally changing phenotypes.  Given the complex nature of traits, a subset of cultivars that show similar phenotypic behavior during one part of the season may possibly diverge at other times.  Secondly, even within a single cultivar, there could be  variability observed across the different seasons, as phenotypic plasticity is a well-established phenomenon in plants \cite{schlichting1986evolution}.
Consequently, it is important for any downstream machine learning workflow to incorporate these complex structural relationships among cultivars in order to improve the efficacy of the prediction tasks. 

In this paper, we model the problems of extracting cross-cultivar  relationships and intra-cultivar variability patterns as one of a structure discovery process. More specifically, given a multi-dimensional point cloud, where each point represents a phenotypic observation of a cultivar in time, we model the problem of identifying structural patterns in data using topological data analysis. In particular, we explore the use of persistent homology \cite{edelsbrunner_persistent_2013}, an active research area within the field of applied algebraic  topology. This is a branch of mathematics that studies \emph{shapes} of data and spaces using algebraic techniques. Topology works with coordinate-free representation of shapes (simplicial complexes) \cite{munkres2018elements},
which are also more robust to small changes in data or to missing data \cite{lum_extracting_2013}.

There are two major classes of techniques within applied  topology---\texttt{mapper} \cite{singh_topological_2007} and persistent homology \cite{edelsbrunner_persistent_2013}. Here, we explore persistent homology to study point clouds as it is better equipped to elucidate topological features that persist over multiple scales (including temporal scales).
Although persistent homology has been widely used in a number of other application domains, it is yet to be explored in any serious depth for agricultural data sets. 
To the best of our knowledge, it has only been applied to decode plant leaf shapes \cite{zhang2021mfcis}. However, the focus of this paper is different; we aim to identify patterns between cultivars and within cultivars based on phenotypic behavior. 

As a concrete application case study, we study the cold hardiness of multiple grape cultivars.
Cold hardiness is a trait that measures how resilient a variety is to cold temperatures. Specialty crops such as grapes and apples can incur a significant loss when the air temperature drops below certain cold hardiness thresholds \cite{mills_cold-hardiness_2006}. However, these thresholds are not fixed, change by the time of the season, and vary by the cultivars. 
Furthermore, due to a large number of cultivars and their divergent behavior across different growing conditions, it becomes important to study: a) the relationships between different cultivars by their cold hardiness trait, and b) any trait level variability as seen in the same cultivar across different seasons. 
Elucidating these relationships will help field scientists devise better frost/cold mitigation protocols customized and effective when applied to different cultivars. 
They could also help us improve the precision of current state-of-the-art cold hardiness prediction models (e.g., \cite{ferguson_modeling_2014}).




\section{Cold Hardiness Data}
\label{sec:data}

This study used the cold hardiness of endo\textendash and ecodormant primary buds from about 30 diverse field-grown grapevine cultivars measured since 1988 at the WSU Irrigated Agriculture Research and Extension Center (IAREC). Locations include the vineyards at IAREC, Prosser, WA (46.29°N lat., -119.74°W long.), the WSU-Roza Research Farm, Prosser, WA (46.25°N lat., -119.73°W long.), and the Ste. Michelle Wine Estates, Paterson, WA (45.96°N lat.; -119.61°W long.). Cane samples containing dormant buds were collected daily, weekly, or at 2-week intervals from leaf fall in October to bud swell in April---i.e., dormant season. The collected samples were analyzed using differential thermal analysis (DTA) \cite{mills_cold-hardiness_2006}. DTA requires putting the samples in a thermoelectric module that senses low-temperature exotherms (LTEs) resulting from the freezing events of individual buds. The module is placed in a controlled chamber, where LTEs  are monitored and registered as the temperature decreases. The result is a measurement of the lethal temperatures at which 10\%, 50\%, and 90\% of the bud population die, denoted by $LTE_{10}$, $LTE_{50}$, and $LTE_{90}$, respectively. Additionally, daily environmental data (e.g., max and min air temperature) from the closest on-site weather station to each vineyard was obtained using the API provided by AgWeatherNet \cite{AgWeatherNet}. 

Thus, for each cultivar, there is a temporal dataset with a varying number of seasons containing daily weather data along with cold hardiness LTE values for the days that samples were collected. For the purpose of this study, a \emph{season} is said to span from September 7th to May 15th---a conservative interval containing the full dormancy period \cite{ferguson_dynamic_2011, ferguson_modeling_2014}. 

\paragraph{Data Summary.} Table \ref{tab:data-description} shows a summary of the number of years of data collected for the five different grape cultivars that have been selected for this study based on their market importance ($n=3,629$ samples). 
Each cultivar dataset contains a row for each day of data collection. 
The key (selected) temporal fields include:
\begin{itemize} 
\item DATE: The date of  observation. 
\item SEASON\_JDAY: Julian day is an integer that represents the count of the number of days since the beginning of the year. For the dormant season, the JDAY continues at the new year---i.e., starts on JDAY 250 (Sep 7th) and ends with JDAY 500 (May 15th). 
\begin{table}[th!]
\centering
\resizebox{1\columnwidth}{!}{
\begin{tabular}{ |l|l|c|c| }
\hline
\multicolumn{1}{|c|}{\textbf{Cultivar}}             &\multicolumn{1}{|c|}{\textbf{LTE Data Seasons}}                      &\makecell{\textbf{Years of} \\ \textbf{LTE Data}}    &\makecell{\textbf{LTE Total} \\ \textbf{Samples}} \\\hline
Cabernet Sauvignon (CS)	 &1988-2022	                                    &34                         &829 \\\hline
Chardonnay (CH)	         &1996-2022	                                    &26                         &783 \\\hline
Concord (CD)	             &1988-'90,'92-'93,'97-'98,'99-2022	&27                         &484 \\\hline
Merlot (MR)	             &1996-2022	                                    &26                         &897 \\\hline
Riesling (WR)             &1988-2022	                                    &34                         &636 \\\hline
\end{tabular}
}
\caption{Summary of LTE data for selected grape cultivars.}
\label{tab:data-description}
\end{table}
\item LTE values at $LTE_{10}$, $LTE_{50}$, and $LTE_{90}$: in $^{\circ}$C. 
\item MIN\_AT, AVG\_AT, MAX\_AT: Minimum, average, and maximum air temperatures respectively, as observed at 1.5 meters above the ground (in degrees Celsius). 
\end{itemize}

\section{Methods}
\label{sec:methods}

\subsection{Persistent Homology}
\label{sec:persistence}

The input to the persistent homology (PH) pipeline \cite{edelsbrunner_persistent_2013} is a high dimensional point cloud $\mathcal{X}$ of $n$ points in $d$ dimensions with a distance $\mathrm{dist}(x,y)$ specified for any pair $x,y \in \mathcal{X}$.
In short, PH characterizes the structure of $\mathcal{X}$ by identifying features in each dimension that persist across multiple scales of $\mathrm{dist}$ values.
It tracks a simplicial complex $K$ built on $\mathcal{X}$ as a function of distance values $\mathrm{dist} \leq r$ for $r \geq 0$.
At any given cutoff $r$, an edge $xy \in K$ when $\mathrm{dist}(x,y) \leq r$.
Similarly, a triangle $xyz \in K$ when every pairwise $\mathrm{dist}$ for points $\{x,y,z\}$ is $\leq r$, and so on.
PH then tracks the evolution of algebraic objects (groups) defined on $K$ as $r$ increases.
It creates a \emph{persistence diagram} (PD) $\mathrm{dgm}_i$ in dimension $i$ that represents each $i$-dimensional feature by a point in the 2D plane with coordinates $\langle$birth,death$\rangle$ corresponding to the values of $r$ at which the feature starts and one at which it stops existing.
In particular, $\mathrm{dgm}_0$ captures the evolution of connected components, while $\mathrm{dgm}_1$ captures the evolution of \emph{holes} in $\mathcal{X}$.
In this work, we concentrate on $\mathrm{dgm}_1$ PDs as holes in $\mathcal{X}$ can be used to capture branching behavior (see Section~\ref{sec:PDforCH}).  

Furthermore, the PH pipeline provides a natural way to compare pairs of data sets $\mathcal{X}$ and $\mathcal{Y}$ based on their branching behavior. 
We first generate the corresponding PDs $\mathrm{dgm}_1(\mathcal{X})$ and $\mathrm{dgm}_1(\mathcal{Y})$.
Considering these PDs as histograms (or probability measures), we compute the \emph{Wasserstein distance} $\mathrm{WD}(\mathrm{dgm}_1(\mathcal{X}), \mathrm{dgm}_1(\mathcal{Y}))$ between them (also known as the Earth mover's distance).
We then use the WD values to directly compare the data sets $\mathcal{X}$ and $\mathcal{Y}$.

\begin{figure*}[tbh!]
\centering
\includegraphics[scale=0.25]{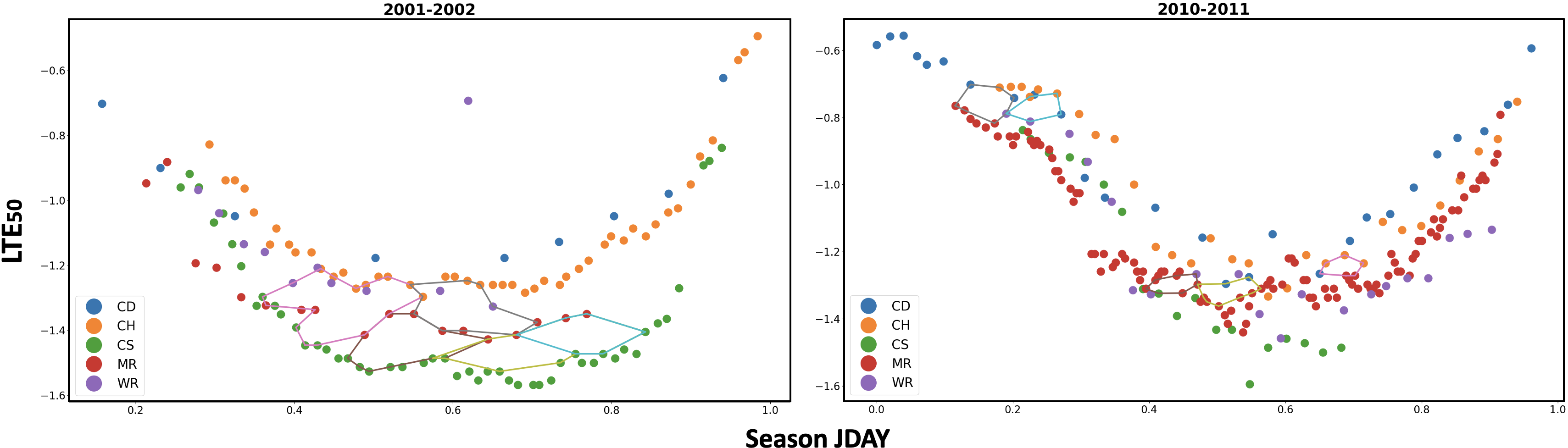} 
\caption{Branching events detected by persistent homology for each point cloud data set. 2001-2002 and 2010-2011 were the pair with the greatest Wasserstein Distance for their respective persistence diagrams among all season pairs.
The y-axis LTE value is calculated using Eqn.~(\ref{eqn:norm}).
}
\label{fig:holesseason}
\end{figure*}

\subsection{Building Persistence Diagrams for Analyzing Cold Hardiness Behavior}
\label{sec:PDforCH}

In the case of cold hardiness, each point in the point cloud $\mathcal{X}$  is a 4-tuple 
$\langle c, s, d, h \rangle$, where $c$ is the cultivar label, $s$ is the season/year, $d$ is a JDAY of the season, and $h$ is the cold hardiness value (i.e., the phenotype) observed that day for that cultivar, which is one of $LTE_{10}$, $LTE_{50}$, or $LTE_{90}$.

Using this input point cloud $\mathcal{X}$, we construct different types of persistence diagrams (Section~\ref{sec:persistence}) to answer different kinds of queries as described below. 

\paragraph{Task 1. Computing Inter-cultivar and Intra-cultivar relationships.}
Consider two cultivars $c_1$ and $c_2$ that exhibit similar LTE values during most of the season except for an interval when they diverge in their values. It is also possible that such divergent behavior may be observable at  multiple time scales---e.g., two cultivars could diverge for days, while another two cultivars could diverge for weeks to months. 

Furthermore, it is  possible the LTE behavior exhibited by a single cultivar across two different seasons is divergent. For instance, a cultivar $c$ may exhibit a higher range of LTE values during one season and a lower range in another season, both during the same JDAY time interval. 

In what follows, we present an approach using persistence diagrams to detect these distinct kinds of divergent behavior.  

\begin{compactenum}[S1)]
\item 
For each cultivar $c$, construct a point cloud $\mathcal{X}_c\subseteq \mathcal{X}$ with all points of the form $\langle d,h \rangle$ taken from all points in $\mathcal{X}$ corresponding to cultivar $c$.
\item
Next, for each cultivar $c$ and using $\mathcal{X}_c$, build a persistence diagram containing $\mathrm{dgm}_0$ (connected components) and $\mathrm{dgm}_1$ (holes) as described in Section~\ref{sec:persistence}. We denote those diagrams for cultivar $c$ as $\mathrm{dgm}_0^c$ and $\mathrm{dgm}_1^c$, respectively.
\item 
We then compare each pair of cultivars by computing the pairwise Wasserstein distance (described in Section~\ref{sec:persistence}) between their respective diagrams.
More specifically, for a cultivar pair $c$ and $c'$, we compute 
$\mathrm{WD}(\mathrm{dgm}_1^{c}(\mathcal{X}_{c}),\mathrm{dgm}_1^{c'}(\mathcal{X}_{c'}))$.
\end{compactenum}

The $\mathrm{dgm}_1$ outputs of step S2 can be used to infer intra-cultivar variability, and the outputs of step S3 to infer inter-cultivar relationships.
Consider a hole detected as part of a $\mathrm{dgm}_1^c(\mathcal{X}_c)$, and let that hole span from day $i$ to day $j$ along the JDAY dimension of the point cloud.
This is indicative of a branching event that starts around day $i$ and ends around day $j$. If the two branching paths are comprised of points from two different seasons $s_1$ and $s_2$ (to be expected), then we infer that the cultivar $c$ shows variable LTE behavior between these two seasons for the JDAY interval $[i,j]$.

Note that different $\mathrm{dist}(.)$ can be used for step S2 to compute the persistence diagrams.
In this paper, we first used two different scaling functions for the two dimensions of the point cloud (JDAY, LTE50), and then used the L2 (Euclidean) distance as the function on the scaled points to construct the persistence diagrams.

For JDAY, we used a simple normalization to scale all points in the range of [0,1]. 

For LTE, instead of using their values directly, we modeled the phenotype value by taking the difference between the observed LTE value and the minimum air temperature recorded on that day (we denote this difference $\delta(s,d,h)$ on day $d$ of season $s$ with LTE $h$). Intuitively, this difference is a strong indicator of the degree of risk that the cultivar faces on that day---as the difference shrinks, the cultivar is at a higher risk.
Note that $\delta(s,d,h)$ can be positive or negative (rare).
Subsequently, we normalize the $\delta$ values as:
\begin{equation}\label{eqn:norm}
\hspace*{-0.05in}
\overline\delta(s,d,h) \hspace*{-0.02in} = \hspace*{-0.02in}
\frac{\delta(s,d,h)-\min_d\{|\delta(s,d,h)|\}\,}{\max_d\{|\delta(s,d,h)|\}-\min_d\{|\delta(s,d,h)|\}} 
\end{equation}
Intuitively, this normalization function is aimed at making an oval or oblong branching shape into a more circular shape, making it suited for hole detection (examples shown in Figure~\ref{fig:holesseason}).


\paragraph{Task 2. Computing seasonal correlations.}
Given multiple seasons, we are also interested in computing pairwise seasonal correlations. Two seasons are said to be \emph{similar} if all the cultivars considered show consistent relative behavior between the two seasons. 
There are two ways to compute this relationship. 
We can directly compare the point clouds for those two seasons. 
Alternatively, we can build the persistence diagrams for those two seasons and compare them. 
We choose the latter approach since persistence diagrams are more compact representations with robust properties \cite{edelsbrunner_persistent_2013}.

\begin{figure*}[tbh!]
\centering
\includegraphics[width=2.1\columnwidth]{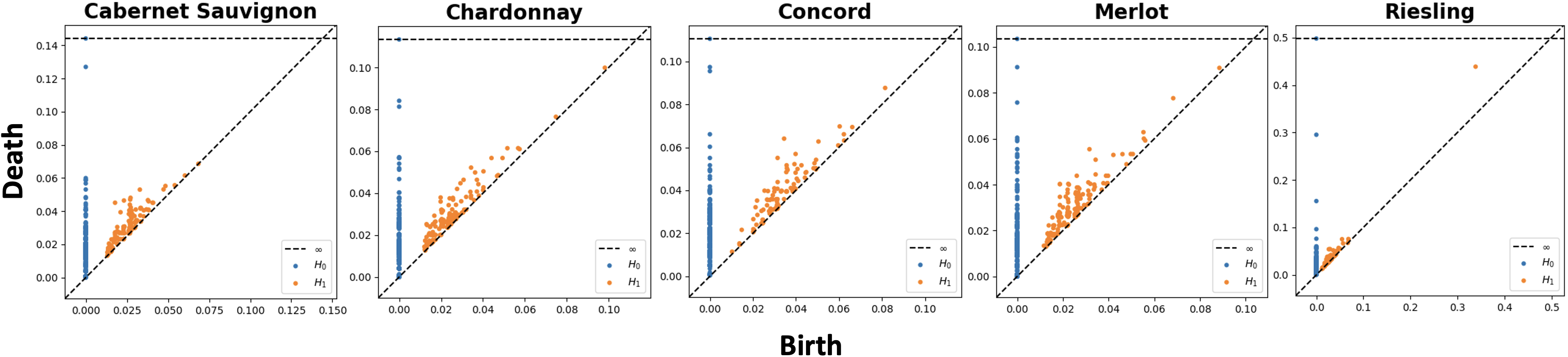} 
\caption{Persistence diagrams for the five grape cultivars, each constructed using their data from seasons 1999 through 2022.
Blue and orange dots represent connected components ($\mathrm{dgm}_0$) and holes ($\mathrm{dgm}_1$) respectively. The vertical distance of a dot from the main diagonal corresponds to the feature's duration. 
}
\label{fig:cultivarspd}
\end{figure*}

\begin{figure*}[tbh]
\centering
\includegraphics[width=2.1\columnwidth]{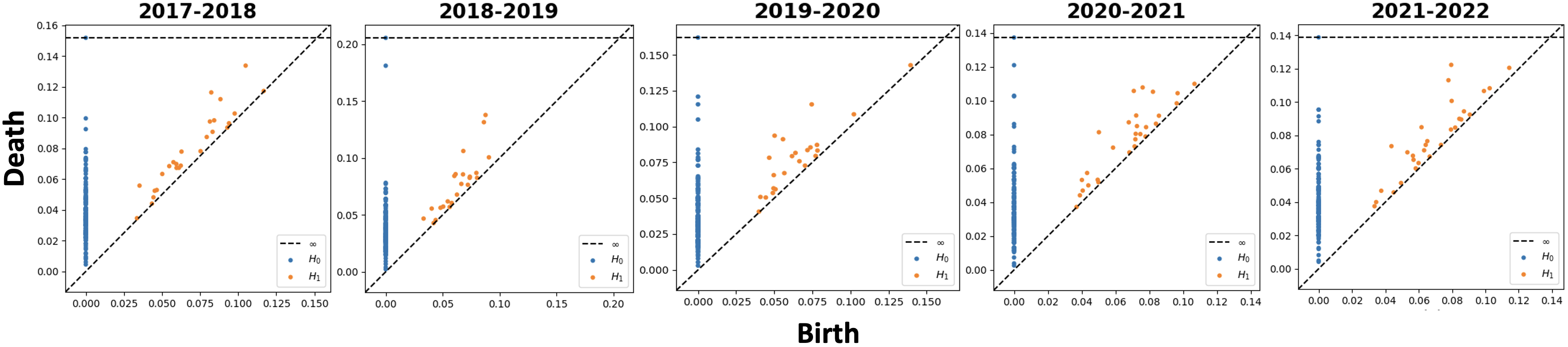} 
\caption{Persistence diagrams for the last 5 seasons of the cold hardiness data. Each diagram is for one season, constructed using all 5 cultivar data for that season.
}
\label{fig:PDpairseasons}
\end{figure*}

\section{Results}
\label{sec:results}

In this section, we present the results of applying our methodology described under Tasks 1 and 2 in Section~\ref{sec:PDforCH} on the grape cold hardiness data set described in Section~\ref{sec:data}.
All experiments shown are for  $LTE_{50}$ values (due to space restrictions). 
All persistence diagrams and Wasserstein distances were computed using the \texttt{Scikit-TDA} package \cite{scikittda2019}. 

\paragraph{Results of Inter-cultivar Comparisons:}
\label{sec:results:intercultivar}
Figure~\ref{fig:cultivarspd} shows the persistence diagrams computed for each cultivar, and
Table~\ref{tab:distance_matrix_cultivars} shows the Wasserstein distance matrix for all cultivar pairs using their persistence diagrams.
As can be seen, CD and MR display the largest distance, while CH and CS constitute the closest pair.
WR (Riesling) has a larger distance to all other cultivars except to CD. 
Note that there are no known connections between the cold hardiness trait and the consumptive type of grape (i.e., wine or juice or table).

\begin{table}[htb]
\centering
\resizebox{0.7\columnwidth}{!}{
\begin{tabular}{ |l|r|r|r|r|r|}
\hline
&	\multicolumn{1}{|c|}{\textbf{CD}}	&	\multicolumn{1}{|c|}{\textbf{CH}}	&	\multicolumn{1}{|c|}{\textbf{CS}}	&	\multicolumn{1}{|c|}{\textbf{MR}}	&	\multicolumn{1}{|c|}{\textbf{WR}}\\\hline
\textbf{CD}	&	0	&	0.236	&	0.201	&	\emph{0.305}	&	0.168\\\hline
\textbf{CH}	&	0.236	&	0	&	\emph{0.146}	&	0.164	&	0.233\\\hline
\textbf{CS}	&	0.201	&	0.146	&	0	&	0.168	&	0.206\\\hline
\textbf{MR}	&	0.305	&	0.164	&	0.168	&	0	&	0.286\\\hline
\textbf{WR}	&	0.168	&	0.233	&	0.206	&	0.286 &	0\\\hline
\end{tabular}
}
\caption{Pairwise distance matrix for cultivars. For each cultivar data from 1999-2022 was used. Each value shows the Wasserstein distance between the $\mathrm{dgm}_1$ obtained for the corresponding two cultivars.
}
\label{tab:distance_matrix_cultivars}
\end{table}

\vspace*{-0.2in}
\paragraph{Results of Intra-cultivar Variability:}
\label{sec:results:intracultivar}
Next, we ask how variable is each cultivar across the different seasons. This is captured by the level of branching observed within the $\mathrm{dgm}_1$ for that cultivar (Task 1, Section~\ref{sec:PDforCH}).
Figure~\ref{fig:cultivarspd}
shows all five persistence diagrams. As can be readily observed, each cultivar has a different profile. 
Intuitively, if a cultivar behaves highly variable from season to season, we can expect to see more branching. However, if those branching events are relatively short-lived, then they correspond to small time scale variations. Longer lived branching events correspond to more persistent divergent behavior. 
From Figure~\ref{fig:cultivarspd}, it can be seen that all three of CH, CD, and MR show numerous holes of wide ranging durations. In contrast, CS and WR show fewer holes with smaller duration. This suggests that CS and WR are less variable compared to the other varieties.

\paragraph{Seasonal comparisons:}
\label{sec:results:seasonalcompare}
We also compared the different seasons (from 1999 to 2022) using the methodology described in Task 2 of Section~\ref{sec:PDforCH}.
Figure~\ref{fig:PDpairseasons} shows the persistence diagrams for only the last 5 seasons (due to space constraints).
We observed that the two most different seasons (i.e., with the largest Wasserstein distance) were the seasons 2001-2002 vs. 2010-2011 (shown 
in Figure~\ref{fig:holesseason}).

\section{Conclusion}
\label{sec:conclusions}
Topological data analysis can be an effective tool to mine for higher order structural information from point cloud data. In this paper, we presented a persistent homology based framework to analyze and glean various types of structural information from a cold hardiness data set. The framework itself is generic and can be extended to other applications within agriculture or other domains.
Future research directions include (but are not limited to) a) using the information gained to improve the prediction accuracy of cold hardiness models; b) adverse testing against noise and incomplete data; and c) exploring ways to use relationships inferred toward data imputation and multi-task learning among cultivars.

\section*{Acknowledgement}
This research was supported by USDA NIFA award No. 2021-67021-35344 (AgAID AI Institute). 
The authors thank the Keller lab at WSU IAREC for data  collection.
\bibliography{aaai23.bib}

\end{document}